\documentclass[12pt,a4paper,twoside]{amsart}
\usepackage{amscd}
\usepackage{amssymb,bm}
\usepackage{undertilde}

\theoremstyle{definition}

\theoremstyle{remark}

\usepackage{amsfonts}

\usepackage{amsmath}
\usepackage{amssymb}
\usepackage{amsfonts}
\usepackage{enumerate}

\numberwithin{equation}{section}

\date{}

\def\PGL{\text{\rm PGL}}
\def\SL{\text{\rm SL}}
\def\GL{\text{\rm GL}}

\def\Ga{\Gamma}
\def\bbr{\mathbb{F}}

\def\bbz{\mathbb{Z}}

\theoremstyle{plain}

\newtheorem*{thm*}{Theorem}

\newtheorem*{prop*}{Proposition}
\newtheorem*{prop**}{\ }
\def\beq{\begin{equation}}

  \def\ee{\end{equation}}

\theoremstyle{definition} 

\newtheorem*{definition*}{Definition}
\newtheorem*{thm1*}{Theorem A1}
\newtheorem*{thm2*}{Theorem A2}

\newtheorem*{conjecture*}{Conjecture}

\newtheorem*{claim*}{Claim}

\newtheorem*{remark*}{Remark}
\def\bbz{\mathbb{Z}}
\def\bbq{\mathbb{Q}}

\def\bbr{\mathbb{R}}
\def\bba{\mathbb{A}}

\def\bbh{\mathbb{H}}
\def\be{\begin{equation}}
\def\ee{\end{equation}}

\def\vare{\varepsilon}

\theoremstyle{remark}  
\begin{document}

\title[Ramanujan Graphs]{Ramanujan Graphs}
\author[A. Lubotzky]{Alexander Lubotzky}
\maketitle


\baselineskip 16pt




Let $X$ be a finite connected $k$-regular graph, $k\ge 3$, with $n$ vertices, and $A$ its adjacency $n \times n$ matrix.  Being symmetric, all its eigenvalues $\lambda$ are real  and it is easy to see that $|\lambda| \le k$, $k$ is always an eigenvalue, and $-k$ is an eigenvalue if and only if $X$ is bi-partite.  The graph $X$ is called \emph{Ramanujan graph} if every eigenvalue $\lambda$ satisfies either $|\lambda| = k$ or $|\lambda| \le 2 \sqrt{k-1}$.  The bound $2\sqrt{k-1}$ is significant: by Alon-Boppana Theorem, (cf. \cite{LPS}) this is the best possible bound one can hope for, for an infinite family of $k$-regular graphs.  The real reason behind  it is as follows:
The universal cover of $X$ (in the sense of algebraic topology) is $\tilde X = T_k$ - the infinite $k$-regular tree.  An old result of Kesten asserts that the spectrum of the adjacency  operator acting on $L^2(T_k)$ is the interval $[-2\sqrt{k-1}, 2\sqrt{k-1}]$.  So, being Ramanujan means for $X$, that all its non-trivial eigenvalues are in the spectrum of its universal cover $\tilde X$.

Ramanujan graphs are optimal expanders from spectral point of view.   Recall that a finite $k$-regular graph $X$ is called $\vare$-expander if  $h(X) \ge \vare$ when $h(X)$ is the Cheeger constant of $X$, namely
$$ h(X) = \min\left\{ |E(Y, \bar Y)|/|Y|\, \Big| {Y \subseteq X,\atop |Y| \le \frac{|X|}{2}}\right\} $$
when $E(Y, \bar Y)$ is the set of edges between $Y$ and its complement.

 Now if we denote $\lambda_1 (X) =\max  \{ \lambda \, | \, \lambda \neq k, \lambda \hbox{ e.v. \ of\ } X\}$, then
 $$ \frac{h^2|X|}{2k} \le \lambda_1 (X) \le 2 h(X) \hbox{\ (cf. \cite{L1})}.$$
 So, Ramanujan graphs are expanders.  Expander graphs are of great importance in combinatorics and computer science (cf. \cite{HLW}
 and the references therein) and also in pure mathematics (cf. \cite{L2}).  Expander graphs serve as  basic building blocks in various network constructions, in many algorithms and so on.  The bound on  their eigenvalues  ensures that the random walk on such graphs converges quickly to the uniform distribution and on Ramanujan graphs this happens in the fastest possible way.  This is one more reason that makes them so useful.

 The existence of Ramanujan graphs is by no means a trivial issue:  While it is known that random $k$-regular graphs are expanders, it is not known if they are Ramanujan.  First examples of infinite families of such graphs were given by explicit construction in \cite{LPS} and \cite{M} for $k = q + 1$, $q$ prime.  In \cite{MSV}, it is shown, by a non constructive method, that for every $k \ge 3$ there exist infinitely many $k$-regular bi-partite Ramanujan graphs.

 Why are Ramanujan graphs named after Ramanujan?  As far as we know Ramanujan had no special interest in graph theory.  Let us explain the reason for this name which was coined in \cite{LPS}.

 Let us look at the following power series
 $$ \triangle (q) = q \prod\limits_{n\ge 1} (1 - q^n) = \Sigma \tau(n) q^n = q - 24q + 252q^3 + \dots $$

 The coefficients $\tau(n)$ define the so called Ramanujan tau function.  Ramanujan conjectured that $\tau(p)\le 2 p^{\frac{11}{2}}$ for every prime $p$.  The importance of $\triangle $ comes from the fact that if we write $q = e^{2\pi i z}$ then $\triangle (z)$ is a cusp form of weight 12 on the upper half plane  ${\bbh} = \{ z = x + iy \, | \, x, y \in \bbr, \; \; y > 0\}$ with respect to the modular group $\Gamma = \SL_2 (\bbz)$ acting on $\bbh$ by Mobius transformation ${ab \choose cd} (z) = \frac{az + b}{cz + d}$. Now if $\Gamma_0 (N) = \{ {ab \choose cd} \in \Gamma \Big| c \equiv 0 (\mod N)\} $ we denote $S_k (N)$ (or more generally $S_k (N, w)$ for a Dirichlet character $w$ of $\bbz/ N\bbz$) the space of cusp forms on $\bbh$ w.r.t. $\Ga_0 (N)$  (and $w$).  The Hecke operators $T_p$ ($p$ prime, $(p, N) = 1$), act, and commute, on each $S_k (N, w)$, and their common eigenfunctions are the Hecke eigenforms.  Now, $S_{12} (\Ga = \Ga_0(1))$ is one dimensional and so $\triangle (z)$ above is such a Hecke eigenform.  Moreover, $\tau(p)$ above is equal to the eigenvalue of $T_p$ acting on $S_{12} (\Gamma)$.  A natural and far reaching generalization of the Ramanujan conjecture mentioned above on the size of $\tau(p)$ is the so called Ramanujan-Peterson (RP) conjecture:  for every Hecke eigenform $f$ in $S_k(N, w)$, the eigenvalues $\lambda_p$ of $T_p$, $ (p, N) = 1$, satisfy $|\lambda_p| \le 2 p^{\frac{k-1}{2}}$.  The reader is referred to \cite{R}   for a concise and clear explanation of all these notions.

The modern approach to automorphic functions via representation theory brought in another point of view on the Ramanujan-Peterson Conjecture.
Satake \cite{S} showed that the RP conjecture is equivalent to the assertion:  Let $\bba = \bba_\bbq$ be the ring of ade$'$les of $\bbq$, and $\pi$ an irreducible cuspidal $\GL_2$-representation in $L^2(\GL_2(\bba)/\GL_2(\bbq))$, such that its component at infinity $\pi_\infty $ is square integrable, then for every prime $p$ the local factor at the $p$-component  $\pi_p$ is a tempered representation.  See \cite{R} for exact definitions.  Here we only mention that a representation  of a (simple) $p$-adic real Lie group $G$ is tempered if it is weakly contained in $L^2(G)$.  The RP conjecture was proved by Deligne (for the special representations that  are relevant to the Ramanujan graphs, the RP  conjecture was actually proven earlier by Eichler).   The representation theoretic formulation suggests  vast generalizations to other algebraic groups.

Let us look at the simple $p$-adic group $G = \PGL_2(\bbq_p)$.  The Bruhat-Tits building associated with $G$ is, in this special case, the $(p+1)$-regular tree $T = T_{p+1}$ which can be identified as $T = G/K$ when $K$ is a maximal compact subgroup of $G$.  If $\Gamma $ is a discrete cocompact subgroup of $G$, then $X = \Gamma \setminus T = \Gamma\setminus G/K$ is a finite $(p+1)$-regular graph.  One can show (see \cite{L1}) that $X$ is a Ramanujan graph if and only if every infinite dimensional $K$-spherical $G$-sub-representation of $L^2(\Gamma\setminus G)$ is tempered.  Deligne theorem, combined with the so called Jacquet-Langlands correspondence, enables the construction of such arithmetic subgroups $\Ga$ for which the temperedness conditions is satisfied and hence Ramanujan graphs are obtained.  This was the method of \cite{LPS} and \cite{M}.  Let us mention that for every $k$, if $G$ is the full automorphism group of $T_k$ and $\Gamma$ a discrete cocompact subgroup of $G$, then $X = \Ga\setminus T_k$ is $k$-regular Ramanujan graph if and only if the same temperedness condition is satisfied:  in other words every non-trivial eigenvalue of $X = \Ga\setminus T_k$ is coming  from the spectrum of $T_k$ if and only if every non-trivial spherical subrepresentation of $L^2(\Ga\setminus G)$ is coming from $L^2(G)$.  This illustrates the connection between the notion of Ramanujan graph and the Ramanujan conjecture.

As mentioned above, the Ramanujan-Peterson conjecture was generalized to other groups, some of its generalizations to $\GL_d$ (instead of only $\GL_2$) led to higher dimensional versions of Ramanujan graphs, the so called Ramanujan complexes.

Finally, another interesting hint to a connection with number theory.  Ihara defined the notion of Zeta function of a $k$-regular graph $X$ and Sunada observed that $X$ is Ramanujan if and only if this Zeta function satisfies ``the Riemann hypothesis" (see \cite{L1} for details).



\newpage

\begin{thebibliography}{JPSH}




\bibitem[HLW]{HLW}   S. Hoory, N. Linial and A. Wigderson, Expander graphs and their applications. \emph{Bull. Amer. Math. Soc. (N.S.)} {\bf 43} (2006), no. 4, 439--561.

\bibitem[L1]{L1}    A. Lubotzky,  Discrete groups, expanding graphs and invariant measures. With an appendix by Jonathan D. Rogawski. \emph{Progress in Mathematics 125}, Birkhäuser Verlag, Basel, 1994. xii+195 pp. ISBN: 3-7643-5075-X

    \bibitem[L2]{L2} A. Lubotzky, Expander graphs in pure and applied mathematics. \emph{ Bull. Amer. Math. Soc. (N.S.)} {\bf 49} (2012), no.~1, 113--162.

  \bibitem[LPS]{LPS}  A. Lubotzky, R. Phillips and P. Sarnak,  Ramanujan graphs. \emph{Combinatorica} {\bf  8} (1988), no. 3, 261--277.

\bibitem[MSV]{MSV}   A.W. Marcus, D.A. Spielman and N. Srivastava,  Interlacing families I: Bipartite Ramanujan graphs of all degrees. \emph{Ann. of Math (2)} {\bf  182} (2015), no. 1, 307-325.

\bibitem [M]{M}  G.A. Margulis,  Explicit group-theoretic constructions of combinatorial schemes and their applications in the construction of expanders and concentrators. (Russian) \emph{Problemy Peredachi Informatsii} {\bf  24}  (1988), no. 1, 51--60.


\bibitem[R]{R}  J.D. Rogawski,     Modular forms, the Ramanujan conjecture and the Jacquet-Langlands correspondence,  An Appendix to \cite{L1}.

\bibitem[S]{S}   I. Satake,  Spherical functions and Ramanujan conjecture. 1966 Algebraic Groups and Discontinuous Subgroups (Proc. Sympos. Pure Math., Boulder, Colo., 1965) pp.~258--264. \emph{Amer. Math. Soc.}, Providence, R.I.

\end{thebibliography}
\end{document}